\documentclass[10pt]{article}

\usepackage{amsmath}
\usepackage{amssymb}
\usepackage{theorem}
\usepackage{xypic}

\renewcommand{\tilde}{\widetilde}

\newcommand{\zz}{{\mathbb Z}}

\newcommand{\pp}{{\mathbb P}}
\newcommand{\cc}{{\mathbb C}}

\renewcommand{\O}{{\cal O}}

\newcommand{\Chow}{\mathop{\rm Chow}\nolimits}

\newcommand{\red}{{\mathop{\rm red}\nolimits}}
\newcommand{\supp}{\mathop{\rm supp}}

\newcommand{\doublearrowstack}[2]%
                      {{{{\scriptstyle#1}\atop{\textstyle\longrightarrow}}\atop{{\textstyle\longrightarrow}\atop{\scriptstyle#2}}}}
\newcommand{\rightleftarrowstack}[2]%
                      {{{{\scriptstyle#1}\atop{\textstyle\longrightarrow}}\atop{{\textstyle\longleftarrow}\atop{\scriptstyle#2}}}}
\newcommand{\leftrightarrowstack}[2]%
                      {{{{\scriptstyle#1}\atop{\textstyle\longleftarrow}}\atop{{\textstyle\longrightarrow}\atop{\scriptstyle#2}}}}

\newtheorem{thm}{Theorem}[section]
\newtheorem{cor}[thm]{Corollary}
\newtheorem{lem}[thm]{Lemma}
\newtheorem{prop}[thm]{Proposition}

\theorembodyfont{\upshape}
\newtheorem{defn}[thm]{Definition}

\newtheorem{rem}[thm]{Remark}
\newtheorem{rmk}[thm]{Remark}

\newenvironment{pf}{\begin{trivlist}\item[]{\sc Proof.}}%
            {\nolinebreak $\Box$ \end{trivlist}}

            {\nolinebreak $\Box$ \end{trivlist}}

\newcommand{\vir}{{\rm vir}}

\renewcommand{\O}{\mathcal O}

\newcommand{\Hilb}{\mathop{\rm Hilb}\nolimits}

\newcommand{\relative}{\, {,} \, }

\newcommand{\noprint}[1]{}

\allowdisplaybreaks[2]

\author{Kai Behrend and Jim Bryan}
\title{Super-rigid Donaldson-Thomas Invariants }

\begin{document}
\sloppy

\maketitle

\abstract{We solve the part of the Donaldson-Thomas theory of
  Calabi-Yau threefolds which comes from super-rigid rational curves.
  As an application, we prove a version of the conjectural
  Gromov-Witten/Donaldson-Thomas correspondence of \cite{MNOP} for
  contributions from super-rigid rational curves. In particular, we
  prove the full GW/DT correspondence for the quintic threefold in
  degrees one and two.}

\tableofcontents

\eject
\section{Introduction}

Let $Y$ be a smooth complex projective Calabi-Yau threefold.  Let
$I_n(Y,\beta)$ be the moduli space of ideal sheaves $I_{Z}\subset \O
_{Y}$, where the associated subscheme $Z$ has maximal dimension equal
to one, the holomorphic Euler characteristic $\chi(\O_Z)$ is equal to
$n$, and the associated 1-cycle has class $\beta\in H_2(Y)$.


Recall that $I_n(Y,\beta)$ has a natural symmetric obstruction theory
\cite{RT}, \cite{BF2}.  
Hence we have the (degree zero) virtual fundamental class
of $I_n(Y,\beta)$, whose degree $N_n(Y,\beta)\in\zz$ is the associated
Donaldson-Thomas invariant.

Let 
\[
C=\sum_{i=1}^s{d_i C_i}
\]
be an effective cycle on $Y$, and assume that the $C_{i}$ are pairwise
disjoint, smoothly embedded rational curves with normal bundle
$N_{C_{i}/Y}\cong \O (-1)\oplus \O (-1)$.  Such curves are called
\emph{super-rigid rational curves} in $Y$
\cite{pa-degenerate,br-pa-gt}. Assume that the class of $C$ is
$\beta$.

Let $J_n(Y,C)\subset I_n(Y,\beta)$ be the locus corresponding to
subschemes $Z\subset Y$ whose associated cycle under the Hilbert-Chow
morphism is equal to $C$ (see Definition~\ref{defnJ}). Since
$J_n(Y,C)\subset I_n(Y,\beta)$ is open and closed (see Remark~\ref{rk:
J is open and closed}), we get an induced virtual fundamental class on
$J_{n} (X,C)$ by restriction. We call
\[
N_n(Y,C)=\deg[J_n(Y,C)]^\vir
\]
the \emph{contribution of $C$ }to the Donaldson-Thomas invariant
$N_n(Y,\beta)$.

The goal of this paper is to compute the invariants $N_n(Y,C)$. 

To formulate our results, we
define a series $P_d(q)\in\zz[[q]]$, for all integers $d\geq0$ by
\begin{equation}\label{eqn: defn of Pd(q)}
\prod_{m=1}^\infty(1+q^mv)^m=\sum_{d=0}^\infty P_d(q)v^d\,.
\end{equation}
Moreover, recall the McMahon function
\begin{equation}\label{eqn: defn of M(q)}
M(q)=\prod_{m=1}^\infty\frac{1}{(1-q^m)^m}\,.
\end{equation}
Then we prove (Theorem~\ref{thm: formula for Z (Y,C)}) that
\[
\sum_{n=0}^\infty N_n(Y,C)q^n=M(-q)^{\chi(Y)}\prod_{i=1}^s
(-1)^{d_{i}} P_{d_i}(-q)\,.
\]

Maulik, Nekrasov, Okounkov, and Pandharipande have conjectured a
beautiful correspondence between Gromov-Witten theory and
Donaldson-Thomas theory which we call the GW/DT correspondence.

As an application of the above formula we prove the GW/DT
correspondence for the contributions from super-rigid rational curves
(Theorem~\ref{thm: super-rigid GW/DT}). In particular, we prove the
full degree $\beta $ GW/DT correspondence (Conjecture 3 of
\cite{MNOP}) for any $\beta $ for which it is known that all cycle
representatives are supported on super-rigid rational curves
(Corollary~\ref{cor: GW/DT holds if beta has only super-rigid
curves}). For example, our results yield the GW/DT correspondence for
the quintic threefold in degrees one and two (Corollary~\ref{cor:
GW/DT holds for quintic in degree one and two}). As far as we know,
these are the first instances of the GW/DT conjecture to be proved for
compact Calabi-Yau threefolds.

The \emph{local} GW/DT correspondence for super-rigid rational curves
follows from the results of \cite{MNOP} as a special case of the
correspondence for toric Calabi-Yau threefolds. In contrast to
Gromov-Witten theory, passing from the local invariants of super-rigid
curves to \emph{global} invariants is non-trivial in Donaldson-Thomas
theory, and can be regarded as the main contribution of this paper.

\subsection{Weighted Euler characteristics}

Our main tool will be the weighted Euler characteristics introduced in
\cite{Beh}. Every scheme $X$ has a canonical $\zz$-valued
constructible function $\nu_X$ on it. The weighted Euler characteristic
$\tilde\chi(X)$ of $X$ is defined as
\[
\tilde\chi(X)=\chi(X,\nu_X)=\sum_{n\in\zz}n\,\chi \left(\nu_X ^{-1}
(n) \right)\,. 
\]
More generally, we use relative weighted Euler characteristics
$\tilde\chi(Z\relative X)$ defined as
\[
\tilde\chi(Z\relative X)=\chi(Z,f^{*}\nu _{X})\,,
\]
for any morphism $f:Z\to  X$. Three fundamental properties are

\begin{itemize}
\item [(i)] if $X\to Y$ is \'etale, then $\tilde\chi(Z\relative X)=\tilde\chi(Z\relative Y)$,
\item [(ii)] if $Z=Z_1\sqcup Z_2$ is a disjoint union,
      $\tilde\chi(Z\relative X)=\tilde\chi(Z_1\relative X)+\tilde\chi(Z_2\relative X)$,
\item [(iii)]
      $\tilde\chi(Z_1\relative X_1)\,\tilde\chi(Z_2\relative X_2)=\tilde\chi\big(Z_1\times
      Z_2\relative X_1\times X_2\big)$.
\end{itemize}

The main result of \cite{Beh}, Theorem~4.18, asserts that if $X$ is a
projective scheme with a symmetric obstruction theory on it, then 
$$\deg[X]^\vir=\tilde\chi(X)\,.$$
Thus we can calculate $N_n(Y,C)$ as $\tilde\chi\left(J_n(Y,C) \right)$. 

We will also need the following fact. If $X$ is an affine scheme with
an action of an algebraic torus $T$ and an isolated fixed point $p\in
X$, and $X$ admits a symmetric obstruction theory compatible with the
$T$-action, then
\[
\nu_X(p)=(-1)^{\dim T_{p}X}\,,
\]
where $T_{p}X$ is the Zariski tangent space of $X$ at $p$. 
This is the main technical result of \cite{BF2}, Theorem~3.4.  

Finally, we will use the following result from \cite{BF2}. If $X$ is a
smooth threefold (not necessarily proper), then
\[
\sum _{m=0}^{\infty }\tilde{\chi } (\Hilb ^{m}X) q^{m}= M (-q)^{\chi (X)}\,.
\]
In the case where $X$ is projective and Calabi-Yau, the above proves
Conjecture 1 of \cite{MNOP}.


\eject
\section{The Calculation}

\subsection{The open subscheme $J_n(Y,C)$}

\begin{defn}\label{defnJ}
Let $C_1,\ldots,C_s$ be pairwise distinct, super-rigid rational curves
on $Y$ and let $(d_1,\ldots,d_s)$ be an $s$-tuple of non-negative
integers.  Let $C=\sum _{i} d_{i}C_{i}$ be the associated 1-cycle on
$Y$ and let $\beta$ be the class of $C$ in homology. Define
\[
J_{n} (Y,C)\subset I_{n} (Y,\beta)
\]
to be the open and closed subscheme consisting of
subschemes $Z\subset Y$ whose associated 1-cycle is equal to $C$.
\end{defn}

\begin{rmk}\label{rk: J is open and closed}
To see that $J_n(Y,C)$ is, indeed, open and closed, consider the
Hilbert-Chow morphism, see \cite{kollar}, Chapter~I, Theorem~6.3, which
is a morphism 
$$f:I_n(Y,\beta)^{sn}\longrightarrow \Chow(Y,d)\,,$$ 
where
$\Chow(Y,d)$ is the Chow scheme of 1-dimensional cycles of degree
$d=\deg\beta$ on $Y$. It is a projective scheme. Moreover,
$I_n(Y,\beta)^{sn}$ is the 
semi-normalization of $I_n(Y,\beta)$. The structure morphism
$I_n(Y,\beta)^{sn}\to I_n(Y,\beta)$ is a homeomorphism of underlying Zariski
topological spaces. Therefore the Hilbert-Chow morphism descends to a
continuous map of Zariski topological spaces
$$|f|:|I_n(Y,\beta)|\longrightarrow |\Chow(Y,d)|\,.$$ Because the
$C_i$ are super-rigid, the cycle $C$ corresponds to an isolated point
of $|\Chow(Y,d)|$.  So the preimage of this point
under $|f|$ is open and closed in $|I_n(Y,\beta)|$. The open subscheme
of $I_n(Y,\beta)$ defined by this open subset is $J_n(Y,C)$.
\end{rmk}

\begin{defn}
As $J_n(Y,C)$ is open in $I_n(Y,\beta)$, it has an induced (symmetric)
obstruction theory and hence a virtual fundamental class of degree
zero. Since $J_n(Y,C)$ is closed in $I_n(Y,\beta)$ it is projective,
and so we can consider the degree of the virtual fundamental class
$$N_n(Y,C)=\deg[J_n(Y,C)]^\vir\,,$$
and call it the {\em contribution of $C$ }to the Donaldson-Thomas
invariant $N_n(Y,\beta)$. 
\end{defn}

\subsection{The closed subset $\tilde J_n(Y,C)$}

\begin{defn}
Let $C=\sum_i d_i C_i$ be as above and denote by $\supp C$ the
reduced closed subscheme of $Y$ underlying $C$. Let 
\[
\tilde{J}_{n} (Y,C)\subset J_{n} (Y,C)\subset  I_{n} (Y,\beta)
\] 
be the closed subset consisting of subschemes $Z\subset
Y$ whose underlying closed subset $Z^\red\subset Y$ is contained in
$\supp C$.
\end{defn}

\begin{rmk}
To see that $\tilde J_n(Y,C)$ is closed in $I_n(Y,\beta)$, let
$W_m\subset Y$ be the $m$-th infinitesimal neighborhood of $\supp
C\subset Y$. For any subscheme $Z\subset Y$, with fixed numerical
invariants $n$ and $\beta $, and such that $Z^{\red }\subset \supp C$,
there exists a sufficiently large $m$ so that $Z\subset W_{m}$.  For
such an $m$, consider the Hilbert scheme $I_n(W_m,\beta)$, which is a
closed subscheme of $I_n(Y,\beta)$, as $W_m$ is a closed subscheme of
$Y$. The underlying closed subset of $I_n(Y,\beta)$ is equal to
$\tilde J_n(Y,C)$.
\end{rmk}

\begin{rmk}
Informally speaking, $J_{n} (Y,C)$ parameterizes subschemes whose one
dimensional components are confined to $C$, but may have embedded
points anywhere in $Y$, whereas $\tilde{J}_{n} (Y,C)$ parameterizes
subschemes where both the one dimensional components and the embedded
points are supported on $C$.
\end{rmk}

\subsection{The open Calabi Yau $N=\O(-1)\oplus\O(-1)$}

We consider the open Calabi-Yau $N$, which is the total space of
the vector bundle $\O(-1)\oplus\O(-1)$ on $\pp^1$. We denote by
$C_0\subset N$ the zero section. We consider the Hilbert scheme
$I_n(N,[dC_0])$. 

Let $\overline{N}$ denote $\pp \big(\O \oplus \O (-1)\oplus \O
(-1)\big)$ and let $D_{\infty }=\overline{N}\backslash N$. Since
$3D_{\infty }$ is an anti-canonical divisor of $\overline{N}$, the
corresponding section defines a trivialization of
$K_{N}$. $\overline{N}$ is naturally a toric variety, $D_{\infty }$ is
an invariant divisor, and we let $T_{0}$ be the subtorus whose
elements act trivially on $K_{N}$. Then $T_{0}$ induces a
$T_{0}$-equivariant symmetric obstruction theory on $I_{n}
(N,[dC_0])$, by Proposition~2.4 of \cite{BF2}. Moreover, the $T_{0}$
fixed points in $I_{n} (N,[dC_0])$ are isolated points whose Zariski
tangent spaces have no trivial factors as $T_{0}$ representations (the
proof of Lemma~4.1, Part (a) and (b) in \cite{BF2} is easily adapted
to prove this).

As in \cite{MNOP}, the $T_{0}$ fixed points in $I_{n} (N,[dC_{0}])$
correspond to subschemes which are given by monomial ideals on the
restriction to the two affine charts of $N$. The number of such fixed
points is given by $p (n,d)$ described below.

Let $p(n,d)$ be the number of triples $(\pi_0,\lambda,\pi_\infty)$,
where $\pi_0$ and $\pi_\infty$ are 3-dimensional partitions and
$\lambda$ a 2-dimensional partition.  The 3-dimensional partitions
each have one infinite leg with asymptotics $\lambda$, and no other
infinite legs. Moreover, $d=|\lambda|$ and $n$ is given by
(\cite{MNOP} Lemma 5)
\[
n=|\pi _{0}|+|\pi _{\infty }|+ \sum _{(i,j)\in \lambda } (i+j+1)\,,
\]
where the size of a three dimensional partition with an infinite leg
of shape $\lambda $ along the $z$ axis is defined by
\[
|\pi | = \# \{(i,j,k)\in \zz_{\geq 0} ^{3} :\quad (i,j,k)\in \pi,  (i,j)\not\in \lambda  \}.
\]

\begin{prop}
We have
$$\tilde\chi\big(I_n(N,[dC_0])\big)=(-1)^{n-d}p(n,d)\,.$$
\end{prop}
\begin{pf}
By Corollary~3.5 of \cite{BF2}, we have
\[
\tilde\chi\big(I_n(N,[dC_0])\big)=\sum_p (-1)^{\dim T_{p}}\,,
\]
where the sum is over all $T_0$-fixed points on $I_n(N,[dC_0])$ and
$T_{p}$ is the Zariski tangent space of $I_n(N,[dC_0])$ at $p$. The
parity of $\dim T_{p}$ can be easily deduced from Theorem~2 of
\cite{MNOP} (just as in the proof of Lemma~4.1~(c) in \cite{BF2}). The
result is $n-d$. So all we have to notice is that $p(n,d)$ is the
number of fixed points of $T_0$ on $I_n(N,[dC_0])$.
\end{pf}

\begin{cor}\label{cor17}
We have
\[
\tilde\chi\Big(\tilde
J_n(N,dC_0)\relative I_n(N,[dC_0])\Big)=(-1)^{n-d}p(n,d)\,.
\]
\end{cor}
\begin{pf}
We just have to notice that all $T_0$-fixed points are contained in
$\tilde J_n(N,dC_0)$. 
\end{pf}

\subsection{The box counting function $p (n,d)$}

Counting three dimensional partitions with given asymptotics has been
shown by Okounkov, Reshetikhin, and Vafa \cite{OkReVa} to be
equivalent to the topological vertex formalism which occurs in
Gromov-Witten theory. They give general formulas for the associated
generating functions in terms of $q$ values of Schur functions which
we will use to prove the following Lemma.

\begin{lem}\label{lem: generating fnc for p (n,d)}
The generating function for $p (n,d)$ is given by
\[
\sum_{n=0}^\infty p(n,d)\,q^n= M(q)^2 P_d(q)\, ,
\]
where the power series $P_{d} (q)$ and $M (q)$ are defined in
Equations~(\ref{eqn: defn of Pd(q)}) and (\ref{eqn: defn of M(q)}).
\end{lem}
\begin{pf}
The generating function for the number of 3-dimensional partitions
with one infinite leg of shape $\lambda $ is given by equation~3.21 in
\cite{OkReVa}:
\[
\sum _{\begin{smallmatrix} \text{3d partitions $\pi
$}\\
\text{asymptotic to $\lambda $} \end{smallmatrix}}q^{|\pi |} = M
(q)q^{-\binom{\lambda }{2}-\frac{|\lambda |}{2}}s_{\lambda ^{t}}
(q^{1/2},q^{3/2},q^{5/2},\dotsc )
\]
where $\lambda ^{t}$ is the transpose partition, $\binom{\lambda
}{2}=\sum _{i}\binom{\lambda _{i}}{2}$, $|\lambda |=\sum _{i}\lambda
_{i}$, and
\[
s_{\lambda^{t} } (q^{1/2},q^{3/2},q^{5/2}, \dotsc )
\]
is the Schur function associated to $\lambda ^{t}$ evaluated at
$x_{i}=q^{(2i-1)/2}$. Using the homogeneity of Schur functions and
writing
\[
s_{\lambda ^{t}} (q)=s_{\lambda ^{t}} (1,q,q^{2},\dotsc )
\]
we can rewrite the right hand side of the above equation as
\[
M (q)q^{-\binom{\lambda }{2}}s_{\lambda ^{t}} (q).
\]
Observing that
\[
\sum _{(i,j)\in \lambda } (i+j+1) = |\lambda | + \binom{\lambda }{2}+\binom{\lambda ^{t}}{2},
\]
we get
\[
\sum _{n,d=0}^{\infty }p (n,d)q^{n}v^{d} = M (q)^{2}\sum _{\lambda }s_{\lambda ^{t}} (q)^{2}q^{|\lambda |+\binom{\lambda ^{t}}{2}-\binom{\lambda }{2}}v^{|\lambda |}.
\]
The hook polynomial formula for $s_{\lambda ^{t}} (q)$
(I.3~ex~2~pg~45,\cite{MacDonald}) is
\begin{equation}\label{eqn: hooklength formula for s(q)}
s_{\lambda ^{t}} (q) = q^{\binom{\lambda }{2}}\prod _{x\in \lambda
^{t}} (1-q^{h (x)})^{-1}
\end{equation}
from which one easily sees that
\[
s_{\lambda ^{t}} (q) = q^{\binom{\lambda }{2}-\binom{\lambda^{t} }{2}}
s_{\lambda } (q).
\]
Therefore
\begin{align*}
\sum _{n,d=0}^{\infty }p (n,d)q^{n}v^{d}&=M (q)^{2}\sum _{\lambda }s_{\lambda }
(q)s_{\lambda ^{t}} (q)q^{|\lambda |}v^{|\lambda |}\\
&=M (q)^{2}\sum _{\lambda }s_{\lambda } (q,q^{2},q^{3},\dotsc )s_{\lambda ^{t}} (v,vq,vq^{2},\dotsc )\\
&=M (q)^{2}\prod _{i,j=1}^{\infty } (1+q^{i+j-1}v)
\end{align*}
where the last equality comes from the orthogonality of Schur
functions (I.4~equation~(4.3$)'$  of \cite{MacDonald}). By rearranging this
last sum and taking the $v^{d}$ term, the lemma is proved.
\end{pf}
\begin{rem}\label{remark: Pd(q) is a ratl fnc}
From the proof of the lemma we see that 
\[
P_{d} (q) = q^{d}\sum _{\lambda \vdash d} s_{\lambda } (q)s_{\lambda ^{t}} (q).
\]
From Equation~(\ref{eqn: hooklength formula for s(q)}), it is immediate
that $P_{d} (q)$ is a rational function in $q$. Moreover, using the
formula for total hooklength (pg~11, I.1~ex~2, \cite{MacDonald}), it
is easy to check that $P_{d} (q)$ is invariant under
$q\mapsto 1/q$.
\end{rem}

\subsection{General $Y$}

\begin{lem}
Let $C$ be a super-rigid rational curve on the Calabi-Yau threefold
$Y$. Then
$$\tilde\chi\big(\tilde
J_n(Y,dC)\relative J_n(Y,dC)\big)=(-1)^{n-d}p(n,d)\,,$$
for all $n$, $d$. 
\end{lem}
\begin{pf}
First of all, by Theorem~3.2 of \cite{laufer}, an analytic
neighborhood of $C$ in $Y$ is isomorphic to an analytic neighborhood
of $C_0$ in $N$.  Therefore, by the analytic theory of Hilbert schemes
(or Douady spaces), see \cite{douady}, we obtain an analytic
isomorphism of $\tilde J_n(Y,dC)$ with $\tilde J_n(N,dC_0)$ which
extends to an isomorphism of a tubular neighborhood of $\tilde
J_n(Y,dC)$ in $I_n(Y,[dC])$ with a tubular neighborhood of $\tilde
J_n(N,dC_0)$ in $I_n(N,[dC_0])$.

The formula for $\nu_X(P)$ in terms of a linking number, 
Proposition~4.22 of \cite{Beh}, shows that $\nu_X(P)$ is an
invariant of the underlying analytic structure of a scheme $X$. Thus,
we have
$$\tilde\chi\big(\tilde
J_n(Y,dC)\relative I_n(Y,[dC])\big)=\tilde\chi\big(\tilde
J_n(N,dC_0)\relative I_n(N,[dC_0])\big)\,.$$
Finally, apply Corollary~\ref{cor17}.
\end{pf}

\begin{lem}\label{stratificationlemma}
Let $f:X\to Y$ be an \'etale morphism  of separated schemes of finite
type over $\cc$. Let $Z\subset X$ be a constructible subset.  Assume
that the restriction of $f$ to the closed points of $Z$, $f:Z(\cc)\to
Y(\cc)$, is injective.  Then we have
$$\tilde\chi\big(f(X)\relative Y\big)=\tilde\chi(Z\relative X)\,.$$
We remark that by Chevalley's theorem (EGA IV, Cor.\ 1.8.5), $f(Z)$ is
constructible, so that $\tilde\chi(f(X)\relative Y)$ is defined.
\end{lem}
\begin{pf}
Without loss of generality, $Z\subset X$ is a closed subscheme and so
$Z\to  Y$ is unramified. 

We claim that there exists a decomposition $Y=Y_1\sqcup\ldots\sqcup
Y_n$ into locally closed subsets, such that, putting the reduced
structure on $Y_i$, the induced morphism $Z_i=Z\times _Y Y_i\to Y_i$ is
either an isomorphism, or $Z_i$ is empty. 

In fact, by generic flatness (EGA IV, Cor.\ 6.9.3), we may assume
without loss of generality that $Z\to Y$ is flat, hence \'etale. By
Zariski's Main Theorem (EGA IV, Cor.\ 18.12.13), we may assume that
$Z\to Y$ is finite, hence finite \'etale. Then, by our injectivity
assumption, the degree of $Z\to Y$ is 1 and so $Z\to Y$ is an
isomorphism. 

Once we have this decomposition of $Y$, the lemma follows from
additivity of the Euler characteristic over such decompositions  and
the \'etale invariance of the canonical constructible function $\nu$. 
\end{pf}

Now we consider the case of a curve with several components. Let
\[
C_{\vec d}=\sum_{i=1}^s d_i C_i
\]
be an effective cycle, where the $C_i$ are pairwise disjoint
super-rigid rational curves in $Y$. We assume $d_i>0$, for all
$i=1,\ldots,s$.

For an $(s+1)$-tuple of non-negative integers
$\vec{m}=(m_0,m_1,\ldots,m_s)$, we let $|\vec{m}|=\sum_{i=0}^sm_i$. 
Consider, for $|\vec{m}|=n$  the open subscheme 
\[
U_{\vec{m}}\subset \Hilb^{m_0}(Y)\times\prod_{i=1}^s J_{m_i}(Y,d_i
C_i)\,,
\]
consisting of subschemes $(Z_0,(Z_i))$ with pairwise disjoint support.  

\begin{lem}
Mapping $(Z_0,(Z_i))$ to $Z_0\cup\bigcup_i Z_i$ defines an \'etale
morphism
\[
f:U_{\vec{m}}\longrightarrow J_n(Y,C_{\vec d})\,.
\]
\end{lem}
\begin{pf}
This is straightforward. See also Lemma~4.4 in \cite{BF}.
\end{pf}

Let us write  $\mathring{Y}$ for $Y\setminus \supp C$ and remark  that
\[
Z_{\vec{m}}=\Hilb^{m_0}(\mathring{Y})\times\prod_{i>0}\tilde 
J_{m_i}(Y,d_i C_i)
\]
is contained in $U_{\vec{m}}$.  Moreover, the
restriction $f:Z_{\vec{m}}\to J_n(Y,C_{\vec d})$ is injective on closed
points.  Finally, every closed point of $J_n(Y,C_{\vec d})$ is
contained in $f(Z_{\vec m})$, for a unique $\vec m$, such that $|\vec
m|=n$.

We will apply Lemma~\ref{stratificationlemma} to the 
diagram%
$$\xymatrix{
{\phantom{x}Z_{\vec m}=\Hilb^{m_0}(\mathring{Y})\times\prod_{i>0}\tilde
J_{m_i}(Y,d_i C_i)}\drto^-{\text{injective}}\dto_{\text{closed 
    subset}} & \\ 
U_{\vec{m}}\rto_{\text{\'etale}}^f\dto_{\text{open embedding}}
&J_n(Y,C_{\vec d})\\ 
{\phantom{MMM}\Hilb^{m_0}(Y)\times\prod_{i>0} J_{m_i}(Y,d_i C_i)}}$$

Thus, we may calculate as follows:
\begin{align*}
\tilde\chi&\big(J_n(Y,C_{\vec d})\big)\\
&=\sum_{|\vec{m}|=n}\tilde\chi\big(f(Z_{\vec m})
\relative J_n(Y,C_{\vec d})\big)\\
&=\sum_{|\vec{m}|=n}\tilde\chi\Big(Z_{\vec m}
\relative U_{\vec{m}}\big)\\
&=\sum_{|\vec{m}|=n}\tilde\chi\Big(\Hilb^{m_0}(\mathring{Y})
\times\prod_{i>0}\tilde J_{m_i}(Y,d_i C_i)\relative \Hilb^{m_0}(Y)
\times\prod_{i>0}J_{m_i}(Y,d_i C_i)\Big)\\
&=\sum_{|\vec{m}|=n}
\tilde\chi\Big(\Hilb^{m_0}(\mathring{Y})\relative \Hilb^{m_0}(Y)\Big)  
\prod_{i>0}\tilde\chi\Big(\tilde J_{m_i}(Y,d_i C_i)\relative 
J_{m_i}(Y,d_i C_i)\Big)\\ 
&=\sum_{|\vec{m}|=n}\tilde\chi\big(\Hilb^{m_0}(\mathring{Y})\big)
\prod_{i>0}(-1)^{m_i-d_i}p(m_i,d_i)\,.
\end{align*}

Now we perform the summation:
\begin{align*}
\sum_{n=0}^\infty &
\tilde\chi\big(J_n(Y,C_{\vec d})\big)q^n\\
&= \sum_{n=0}^\infty\Bigg(\sum_{|\vec{m}|=n}
\tilde\chi\big(\Hilb^{m_0}(\mathring{Y})\big)
\prod_{i=1}^s(-1)^{m_i-d_i}p(m_i,d_i)\Bigg)q^n\\
&= \sum_{n=0}^\infty\sum_{|\vec{m}|=n}
\tilde\chi\big(\Hilb^{m_0}(\mathring{Y})\big)q^{m_0}
\prod_{i=1}^s(-1)^{d_i}p(m_i,d_i)(-q)^{m_i}\\
&=\Bigg(\sum_{m_0=0}^\infty
\tilde\chi\big(\Hilb^{m_0}(\mathring{Y})\big)q^{m_0}\Bigg) 
\sum_{\vec{m}'}
\prod_{i=1}^s (-1)^{d_i}p(m_i,d_i)(-q)^{m_i}\\
&=M(-q)^{\chi(\mathring{Y})}
\prod_{i=1}^s (-1)^{d_i}\sum_{m_i=0}^\infty
p(m_i,d_i)(-q)^{m_i}\\
&=M(-q)^{\chi(\mathring{Y})}
\prod_{i=1}^s M(-q)^2 (-1)^{d_{i}} P_{d_i}(-q)\\
&=M(-q)^{\chi(\mathring{Y})}M(-q)^{2s}
\prod_{i=1}^s (-1)^{d_{i}} P_{d_i}(-q)\\
&=M(-q)^{\chi({Y})}
\prod_{i=1}^s (-1)^{d_{i}} P_{d_i}(-q)
\end{align*}

By the main result of \cite{Beh}, Theorem~4.18, we have
$$N_n(Y,C_{\vec d})=\tilde\chi\big( J_n(Y,C_{\vec d})\big)\,.$$
This finishes the proof of:

\begin{thm}\label{thm: formula for Z (Y,C)}
The partition function for the contribution of $C_{\vec d}$ to the
Donaldson-Thomas invariants of $Y$ is given by
\[
Z(Y,C_{\vec d})=\sum_{n=0}^\infty N_n(Y,C_{\vec d})\,q^n=M(-q)^{\chi(Y)}\prod_{i=1}^s
(-1)^{d_{i}} P_{d_i}(-q)\,.
\]
\end{thm}

\begin{cor}\label{cor: formula for Z'(Y,C)}
Define the reduced partition function
$$Z'(Y,C_{\vec d})=\frac{Z(Y,C_{\vec d})}{Z (Y,0)}\,.$$
Then we have
\[
Z'(Y,C_{\vec d})=\prod_{i=1}^s (-1)^{d_{i}} P_{d_i}(-q)\,,
\]
a rational function in $q$, invariant under $q\mapsto
1/q$.
\end{cor}
\begin{pf}
Behrend and Fantechi prove \cite{BF2} that 
\[
Z (Y,0) = \sum _{n=0}^{\infty }\tilde{\chi } (\Hilb ^{n}Y)q^{n}=M (-q)^{\chi (Y)};
\]
the formula for $Z' (Y,C_{\vec d})$ then follows immediately from
Theorem~\ref{thm: formula for Z (Y,C)}. For the proof that $Z'
(Y,C_{\vec d})$ is a rational function invariant under $q\mapsto 1/q$,
see Remark~\ref{remark: Pd(q) is a ratl fnc}.
\end{pf}

\eject
\section{The super-rigid GW/DT
  correspondence.}

\subsection{The usual GW/DT correspondence}
The Gromov-Witten/Donaldson-Thomas correspondence of \cite{MNOP} can
be formulated as follows.

Let $Y$ be a Calabi-Yau threefold and let
\[
Z_{DT} (Y,\beta ) = \sum _{n\in \zz }N_{n} (Y,\beta )q^{n}
\]
be the partition function for the degree $\beta $ Donaldson-Thomas
invariants. Let
\[
Z'_{DT} (Y,\beta )=\frac{Z_{DT} (Y,\beta )}{Z_{DT} (Y,0)}
\]
be the reduced partition function.

In Gromov-Witten theory, the reduced partition function for the degree
$\beta $ Gromov-Witten invariants, $Z'_{GW} (Y,\beta )$, is given by
the coefficients of the exponential of the $\beta \neq 0$ part of the
potential function:
\begin{equation}\label{eqn: def'n of GW potential}
1+\sum _{\beta \neq 0}Z'_{GW} (Y,\beta )v^{\beta } = \operatorname{exp} \left(\sum _{\beta \neq 0}N_{g}^{GW} (Y,\beta )u^{2g-2}v^{\beta } \right).
\end{equation}
Here
\[
N_{g}^{GW} (Y,\beta ) = \operatorname{deg}[\overline{M} _{g} (Y,\beta
)]^{\vir }
\]
is the genus $g$, degree $\beta $ Gromov-Witten invariant of $Y$.

The conjectural GW/DT correspondence states that 
\begin{itemize}
\item [(i)] The degree 0 partition function in Donaldson-Thomas theory
      is given by
\[
Z_{DT} (Y,0) = M (-q)^{\chi (Y)},
\]
\item [(ii)] $Z'_{DT} (Y,\beta)$ is a rational function in $q$, invariant
      under $q\mapsto 1/q$, and
\item [(iii)] the equality
\[
Z'_{GW} (Y,\beta ) = Z'_{DT} (Y,\beta )
\] 
holds under the change of variables $q=-e^{iu}.$
\end{itemize} 

Part (i) is proved for all $Y$ in \cite{BF2}.

\subsection{The super-rigid GW/DT correspondence}
In an entirely parallel manner, we can formulate the GW/DT
correspondence for $N_{n} (Y,C_{\vec d})$, the contribution from a collection
of super-rigid rational curves $C_{\vec d}=\sum _{i}d_{i}C_{i}$. 

Just as in Donaldson-Thomas theory, there is an open component of the
moduli space of stable maps
\[
\overline{M}_{g} (Y,C_{\vec d})\subset \overline{M}_{g} (Y,\beta )
\]
parameterizing maps whose image lies in the support of $C_{\vec
d}$. There are corresponding invariants given by the degree of the
virtual class:
\[
N^{GW}_{g} (Y,C_{\vec d})=\operatorname{deg}[\overline{M}_{g}
(Y,C_{\vec d})]^{\vir}
\]
We define $Z'_{GW} (Y,C_{\vec d})$ by replacing $N_{g}^{GW} (Y,\beta
)$ on the right side of formula~\eqref{eqn: def'n of GW potential} by
$N_{g}^{GW} (Y,C_{\vec d})$.

Then we can formulate our results as follows. 
\begin{thm}\label{thm: super-rigid GW/DT}
The GW/DT correspondence holds for the contributions from super-rigid
rational curves. Namely, let $C_{\vec d}=d_{1}C_{1}+\dotsb
+d_{s}C_{s}$ be a cycle supported on pairwise disjoint super-rigid
rational curves 
$C_{i}$ in a Calabi-Yau threefold $Y$, and let $Z'_{DT} (Y,C_{\vec d})$ and
$Z'_{GW} (Y,C_{\vec d})$ be defined as above. Then
\begin{itemize}
\item [(ii)] $Z'_{DT} (Y,C_{\vec d})$ is a rational function of $q$,
      invariant under $q\mapsto 1/q$, and
\item [(iii)] the equality
\[
Z'_{DT} (Y,C_{\vec d})=Z'_{GW} (Y,C_{\vec d})
\]
holds under the change of variables $q=-e^{iu}$.
\end{itemize}
\end{thm}
\begin{pf} 
For (ii), see Corollary~\ref{cor: formula for Z'(Y,C)}. To prove
(iii), we reproduce a calculation well known to the experts
(e.g. \cite{Katz}).

By the famous multiple cover formula of Faber-Pandharipande
\cite{Fa-Pa} (see also \cite{pa-degenerate}),
\[
N_{g}^{GW} (Y,C_{\vec d}) = \sum _{i=1}^{s} c(g,d_{i})\,,
\]
where $c (g,d)$ is given by 
\[
\sum _{g\geq 0} c (g,d) u^{2g-2} = \frac{1}{d}\left(\sin
\left(\frac{du}{2} \right) \right)^{-2}.
\]
We compute $Z'_{GW} (Y,C_{\vec d}) $  and make the substitution $q=-e^{iu}$:
\begin{align*}
1+\sum _{(d_{1},\dotsc ,d_{s})\neq 0}^{\infty } & Z'_{GW} (Y,C_{\vec d})v_{1}^{d_{1}}\dotsb v_{s}^{d_{s}}  \\
&=\exp \left(\sum _{j=1}^{s} \sum _{d_{j}=1}^{\infty }\sum _{g=0}^{\infty } c (g,d_{j})u^{2g-2}v_{j}^{d_{j}} \right)\\
&=\prod _{j=1}^{s}\exp \left( \sum _{d_{j}=1}^{\infty }\frac{v_{j}^{d_{j}}}{d_{j}}\left(2\sin \frac{d_{j}u}{2} \right)^{-2}\right)\\
&=\prod _{j=1}^{s}\exp \left(\sum _{d_{j}=1}^{\infty
}\frac{v_{j}^{d_{j}}}{d_{j}}\frac{-e^{id_{j}u}}{\left(1-e^{id_{j}u} \right)^{2}} \right)\\
&=\prod _{j=1}^{s}\exp \left(\sum _{d_{j}=1}^{\infty }\sum _{m_{j}=1}^{\infty }\frac{-m_{j}}{d_{j}}e^{id_{j}m_{j}u}v_{j}^{d_{j}} \right)\\
&=\prod _{j=1}^{s}\exp \left(\sum _{m_{j}=1}^{\infty }m_{j}\log\left(1-v_{j}e^{im_{j}u} \right) \right)\\
&=\prod _{j=1}^{s}\prod _{m_{j}=1}^{\infty }\left(1- (-q)^{m_{j}}v_{j} \right)^{m_{j}}\\
&=\prod _{j=1}^{s}\sum _{d_{j}=0}^{\infty }P_{d_{j}} (-q) (-v_{j})^{d_{j}}\\
&=\sum _{(d_{1},\dotsc ,d_{s})}\prod _{j=1}^{s} (-1)^{d_{j}}P_{d_{j}}
(-q)v_{j}^{d_{j}}.
\end{align*}
Therefore, 
\[
Z'_{GW} (Y,C_{\vec d}) = \prod _{j=1}^{s} (-1)^{d_{j}}P_{d_{j}} (-q)
\]
and so by comparing with Corollary~\ref{cor: formula for Z'(Y,C)} the
theorem is proved.
\end{pf}

The following corollary is immediate.

\begin{cor}\label{cor: GW/DT holds if beta has only super-rigid curves}
Let $Y$ be a Calabi-Yau threefold and let $\beta \in H_{2} (Y,\zz)$ be
a curve class such that all cycle representatives of $\beta $ are
supported on a collection of pairwise disjoint, super-rigid rational
curves. Then the degree $\beta $ GW/DT correspondence holds:
\[
Z'_{DT} (Y,\beta ) = Z'_{GW} (Y,\beta ).
\]
\end{cor}

For example, we have:

\begin{cor}\label{cor: GW/DT holds for quintic in degree one and two}
Let $Y\subset \pp ^{4}$ be a quintic threefold, and let $L$ be
the class of the line. Then for $\beta $ equal to $L$ or $2L$, the
GW/DT correspondence holds.
\end{cor}
\begin{pf}
By deformation invariance of both Donaldson-Thomas and Gromov-Witten
invariants, it suffices to let $Y$ be a generic quintic threefold. It
is well known that there are exactly 2875 pairwise disjoint lines on
$Y$ and they are all super-rigid. The conics on $Y$ are all planar and
hence rational, and it is known that there are exactly 609250 pairwise
disjoint conics and they are super-rigid as well. For these facts and more,
see \cite{Katz86compmath}.
\end{pf}

Note that we cannot prove the GW/DT conjecture by this method for the
quintic in 
degree three (and higher) due to the presence of elliptic curves in
degree three.

Explicit formulas for the reduced Donaldson-Thomas partition
function of a generic quintic threefold $Y$ in degrees one and two are
given below:
\begin{align*}
Z'_{DT} (Y,L)&= 2875\frac{q}{(1-q)^{2}}\\
Z'_{DT} (Y,2L)&= 609250\frac{q}{(1-q)^{2}}\cdot 2875\cdot
\frac{-2q^{3}}{(1+q)^{4} (1-q)^{2}}\\
&=-3503187500\frac{1}{(q-q^{-1})^{4}}
\end{align*}

\eject


\end{document}